\documentclass[12pt]{amsart}
\usepackage{amsmath,amssymb,amsthm,amscd}
\usepackage{epsfig} 

\theoremstyle{remark}

\newcommand{\C}{{\mathbb C}}

\newcommand{\U}{{\mathfrak U}}
\newcommand{\UU}{{\mathbf U}}
\renewcommand{\H}{\mathbf{H}}

\renewcommand{\ge}{\geqslant}
\renewcommand{\le}{\leqslant}
\newcommand{\gl}{\mathfrak{gl}}
\renewcommand{\sl}{\mathfrak{sl}}
\newcommand{\so}{\mathfrak{so}}
\renewcommand{\sp}{\mathfrak{sp}}
\newcommand{\GL}{\mathsf{GL}}
\newcommand{\SL}{\mathsf{SL}}
\newcommand{\Sp}{\mathsf{Sp}}
\renewcommand{\O}{\mathsf{O}}
\newcommand{\SO}{\mathsf{SO}}
\newcommand{\End}{\operatorname{End}}

\newcommand{\g}{\mathfrak{g}}

\newcommand{\Z}{{\mathbb Z}}
\newcommand{\Q}{{\mathbb Q}}
\renewcommand{\S}{\mathfrak{S}}
\renewcommand{\SS}{\mathbf{S}}
\newcommand{\A}{\mathcal{A}}
\newcommand{\B}{\mathcal{B}}

\allowdisplaybreaks

\begin{document}
\title{New versions of Schur-Weyl duality}
\author{Stephen Doty}
\address{Mathematics and Statistics, Loyola University Chicago, 
 Chicago, Illinois 60626 U.S.A.}
\thanks{These notes are based on a lecture, various versions of which I have
given in the past year, in a number of locations, including Stuttgart,
Birmingham, Queen Mary (London), Lancaster, Manchester, Oxford, and
Cambridge. I'm grateful to the organizers of those events for the
opportunity to present these ideas.}
\email{doty@math.luc.edu}
\begin{abstract}
After reviewing classical Schur-Weyl duality, we present some other
contexts which enjoy similar features, relating to Brauer algebras
and classical groups.
\end{abstract}

\maketitle

\section{Classical Schur-Weyl duality}

\subsection{Schur's double-centralizer result}\label{1.1}
Consider the vector space $V=\C^n$. The symmetric group $\S_r$ acts
naturally on its $r$-fold tensor power $V^{\otimes r}$, by permuting
the tensor positions. This action obviously commutes with the natural
action of $\GL_n = \GL_n(\C)$, acting by matrix multiplication in each
tensor position. So we have a $\C\GL_n$-$\C\S_n$ bimodule structure on
$V^{\otimes r}$.  (Here $\C G$ denotes the group algebra of a group
$G$.) In 1927, Schur \cite{S} proved that the image of each group
algebra under its representation equals the full centralizer algebra
for the other action. More precisely, if we name the representations
as follows
\begin{equation}\label{1}
\begin{CD}
\C\GL_n @>\rho>> \End(V^{\otimes r}) @<\sigma<< \C\S_r 
\end{CD}
\end{equation}
then we have equalities
\begin{align}
\rho(\C\GL_n) &= \End_{\S_r}(V^{\otimes r}) \label{a}\\
\sigma(\C\S_r) &= \End_{\GL_n}(V^{\otimes r})\label{b}.
\end{align}
(Here, for a given set $S$ operating on a vector space $T$ through
linear endomorphisms, $\End_S(T)$ denotes the set of linear
endomrphisms of $T$ commuting with each endomorphism coming from $S$.)

Results of Carter-Lusztig \cite{CL} and J.A.\ Green \cite{G} (and
others) show that all the above statements remain true if one replaces
$\C$ by an arbitrary infinite field $K$.

\subsection{Schur algebras}\label{1.1a}
The finite-dimensional algebra in \eqref{a} above, for any $K$, is
known as the {\em Schur algebra}, and often denoted by $S_K(n,r)$ or
simply $S(n,r)$. The Schur algebra ``sees'' the part of the rational
representation theory of the algebraic group $\GL_n(K)$ occurring (in
some appropriate sense) in $V^{\otimes r}$. More precisely, there is
an equivalence between $r$-homogeneous polynomial representations of
$\GL_n(K)$ and $S_K(n,r)$-modules. In characteristic $0$, those
representations (as $r$ varies) determine all finite-dimensional
rational representations, while in positive characteristic they still
provide a tremendous amount of information.

The representation $\sigma$ in \eqref{1} is faithful if $n\ge r$, so
$\sigma$ induces an isomorphism
\begin{equation}
K\S_r \simeq \End_{\GL_n}(V^{\otimes r}) = \End_{S_K(n,r)}(V^{\otimes r})
\qquad (n\ge r).
\end{equation}
This leads to intimate connections between polynomial representations
of $\GL_n(K)$ and representations of $K\S_r$, a theme that has been
exploited by many authors in recent years. Perhaps the most dramatic
example of this is the result of Erdmann \cite{E} (building on
previous work of Donkin \cite{Donkin:Tilt} and Ringel \cite{Ringel})
which shows that knowing decomposition numbers for all symmetric
groups in positive characteristic will determine the decomposition
numbers for general linear groups in the same
characteristic. Conversely, James \cite{J} had already shown that the
decomposition matrix for a symmetric group is a submatrix of the
decomposition matrix for an appropriate Schur algebra.  Thus the
(still open) general problem of determining the modular characters of
symmetric groups is equivalent to the similar problem for general
linear groups (over infinite fields).

\subsection{The enveloping algebra approach} \label{1.2}
Return to the basic setup, over $\C$. One may differentiate the action
of the Lie group $\GL_n(\C)$ to obtain an action of its Lie algebra
$\gl_n$.  Replacing the representation $\rho$ in \eqref{1} by its
derivative representation $d\rho: \U(\gl_n) \to \End(V^{\otimes r})$
leads to the following alternative statement of Schur's result:
\begin{align}
d\rho(\U(\gl_n)) &= \End_{\S_r}(V^{\otimes r}) \label{d}\\
\sigma(\C\S_r) &= \End_{\gl_n}(V^{\otimes r})\label{e}.
\end{align}
In particular, the Schur algebra (over $\C$) is a homomorphic image of
$\U(\gl_n)$.  All of this works over an arbitrary integral domain $K$
if we replace $\U(\gl_n)$ by its ``hyperalgebra'' $\U_K := K\otimes_\Z
\U_\Z$ obtained by change of ring from a suitable $\Z$-form of
$\U(\gl_n)$; see \cite{Donkin:SA1}.  (One can adapt the Kostant
$\Z$-form, originally defined for the enveloping algebra of a {\em
semisimple} Lie algebra, to the {\em reductive} $\gl_n$.)

\subsection{The quantum case}\label{1.3}
Jimbo \cite{Jimbo} extended the results of \ref{1.2} to the quantum
case (where the quantum parameter is not a root of unity). One needs
to replace $\S_r$ by the Iwahori-Hecke algebra $\H(\S_r)$ and replace
$\U(\gl_n)$ by the quantized enveloping algebra $\UU(\gl_n)$. The
analogue of the Schur algebra in this context is known as the
$q$-Schur algebra, often denoted by $\SS(n,r)$ or $\SS_q(n,r)$. Dipper
and James \cite{DJ} have shown that $q$-Schur algebras are fundamental
for the modular representation theory of {\em finite} general linear
groups.

As many authors have observed, the picture in \ref{1.1} can also be
quantized.  For that one needs a suitable quantization of the
coordinate algebra of the algebraic group $\GL_n$.

There is a completely different (geometric) construction of $q$-Schur
algebras given in \cite{BLM}.

\subsection{Integral forms}
The Schur algebras $S_\C(n,r)$ admit an integral form $S_\Z(n,r)$ such
that $S_K(n,r) \simeq K \otimes_\Z S_\Z(n,r)$ for any field $K$. In
fact $S_\Z(n,r)$ is simply the image of $\U_\Z$ (see \ref{1.2}) under
the surjective homomorphism $\U(\gl_n) \to S_\C(n,r)$.  Similarly, the
quantum Schur algebra $\SS_{\Q(v)}(n,r)$ admits an integral form
defining all specializations via base change. One needs to replace $\Z$
by $\A = \Z[v,v^{-1}]$; then the integral form $\SS_\A(n,r)$ is the
image of the Lusztig $\A$-form $\UU_\A$ under the surjection
$\UU(\gl_n) \to \SS_{\Q(v)}(n,r)$. (To match this up with various
specializations in the literature, one often has to take $q = v^2$.)

\subsection{Generators and relations}
Recently, in joint work with Giaquinto (see \cite{DG}), a very simple
set of elements generating the kernel of the surjection $\U(\gl_n) \to
S_\C(n,r)$ was found.  A very similar set of elements generates the
kernel of the surjection $\UU(\gl_n) \to \SS_{\Q(v)}(n,r)$.  These
elements are expressible entirely in terms of the Chevalley generators
for the zero part of $\U(\gl_n)$ or $\UU(\gl_n)$. Thus we obtain a
presentation of $S_\C(n,r)$ and $\SS_{\Q(v)}(n,r)$ by generators and
relations, compatible with the usual Serre (Drinfeld-Jimbo)
presentation of $\U(\gl_n)$ (resp., $\UU(\gl_n)$).  As a result, we
find a certain subset of the integral PBW-basis for $\U(\gl_n)$
or $\UU(\gl_n)$ the image of which gives an integral basis for
$S_\Z(n,r)$ or $\SS_{\A}(n,r)$.  This basis yields a similar basis in
any specialization. Moreover, a subset of it provides a new integral
basis of $\H(\S_n)$.

\section{The Brauer algebra}

From now on I will assume, unless stated otherwise, that the
underlying field is $\C$ (it could just as well be any field of
characteristic zero). One expects that many statements will be valid
over an arbitrary infinite field, via some appropriate integral form,
similar to what happens in type $A$.

\subsection{The algebra $\B_r^{(x)}$} \label{2.1}
\begin{figure} [ht]
       \centering
       \includegraphics{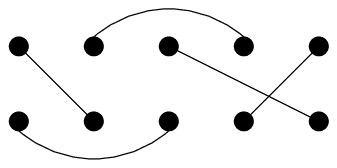}\qquad
       \includegraphics{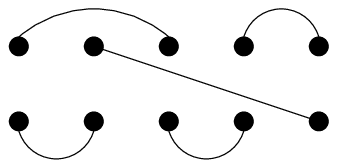}
       \caption{Two Brauer diagrams $D_1$, $D_2$ for $r=5$.}
       \label{fig1}
       \end{figure}

Let $R$ be a commutative ring, and consider the free $R[x]$-module
$\B_r^{(x)}$ with basis consisting of all $r$-diagrams.  An
$r$-diagram is an (undirected) graph on $2r$ vertices and $r$ edges
such that each vertex is incident to precisely one edge.  One usually
thinks of the vertices as arranged in two rows of $r$ each, the {\em
top} and {\em bottom} rows. (See Figure \ref{fig1}.) Edges connecting
two vertices in the same row (different rows) are called {\em
horizontal} (resp., {\em vertical}).  We can compose two such diagrams
$D_1$, $D_2$ by identifying the bottom row of vertices in the first
diagram with the top row of vertices in the second diagram. The result
is a graph with a certain number, $\delta(D_1,D_2)$, of interior
loops. After removing the interior loops and the identified vertices,
retaining the edges and remaining vertices, we obtain a new
$r$-diagram $D_1 \circ D_2$, the composite diagram.

\begin{figure} [ht]
       \centering
       \includegraphics{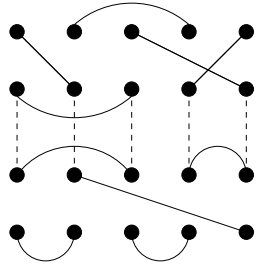} 
       \caption{Computing the composite of $D_1$ and $D_2$.}
       \label{fig2}
       \end{figure}

\begin{figure} [ht]
       \centering
       \includegraphics{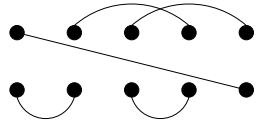}
       \caption{The composite diagram $D_1 \circ D_2$.}
       \label{fig3}
       \end{figure}

Multiplication of $r$-diagrams is defined by the rule
$$
D_1 \cdot D_2 = x^{\delta(D_1,D_2)} (D_1 \circ D_2).
$$
One can check that this multiplication makes $\B_r^{(x)}$ into an
associative algebra; this is the Brauer algebra. (See Figures
\ref{fig1}--\ref{fig3} for an illustration of the multiplication in
the Brauer algebra.)

Note that if we take $x=1$ then the set of $r$-diagrams is a monoid
under diagram composition, and $\B_r^{(1)}$ is simply the semigroup
algebra of that monoid.

For any $x$ the group algebra $R[x]\S_r$ may be identified with the
subalgebra of $\B_r^{(x)}$ spanned by the diagrams containing only
vertical edges. Such Brauer diagrams provide a graphical depiction of
permutations. The group algebra $R[x]\S_r$ of $\S_r$ also appears as a
quotient of $\B_r^{(x)}$, the quotient by the two-sided ideal spanned
by all diagrams containing at least one horizontal edge.

Label the vertices in each row of an $r$-diagram by the indices $1,
\dots, r$. For any $1 \le i \ne j \le r$ let $c_{i,j}$ be the
$r$-diagram with horizontal edges connecting vertices $i$, $j$ on the
top and bottom rows. All other vertices in the diagram $c_{i,j}$ are
vertical, connecting vertex $k$ on the top and bottom rows, for all $k
\ne i, j$. Brauer observed that $\B_r^{(x)}$ is generated by the
permutation diagrams together with just one of the $c_{i,j}$.

\subsection{Schur-Weyl duality}\label{2.2}
Brauer \cite{Br} introduced the algebra $\B_r^{(x)}$ in 1936 to
describe the invariants of symplectic and orthogonal groups acting on
$V^{\otimes r}$. (Brauer's conventions were slightly different; we are
here following the approach of Hanlon and Wales \cite{HW}, who pointed
out that $\B_r^{(-n)}$ is isomorphic with the algebra defined by
Brauer to deal with the symplectic case.)  Let $G$ be $\Sp_n$ or
$\O_n$, where $n$ is even in the first instance. By restricting the
action $\rho$ considered in \ref{1.1} we have an action of $G$ on
$V^{\otimes r}$. One can extend the action of $\S_r$ to an action of
$\B_r^{(\epsilon n)}$ (over $\C$) on $V^{\otimes r}$, where $\epsilon
= -1$ if $G=\Sp_n$ and $\epsilon = 1$ if $G=\O_n$.  To do this, it is
enough to specify the action of the diagram $c_{i,j}$.  This acts on
$V^{\otimes r}$ as one of Weyl's {\em contraction} maps contracting in
tensor positions $i$ and $j$.  So we have (commuting) representations
\begin{equation}\label{7}
\begin{CD}
\C G @>\rho>> \End(V^{\otimes r}) @<\sigma<< \B_r^{(\epsilon n)} 
\end{CD}
\end{equation}
which satisfy Schur-Weyl duality; i.e., the image of each
representation equals the full centralizer algebra of the other
action:
\begin{align}
\rho(\C G) &= \End_{\B_r^{(\epsilon n)}}(V^{\otimes r}) \label{7a}\\
\sigma(\B_r^{(\epsilon n)}) &= \End_{G}(V^{\otimes r}). \label{7b}
\end{align}
The algebras in equality \eqref{7a} are the symplectic and orthogonal
Schur algebras (see \cite{Donkin:Arcata}, \cite{D:PRAMSACT},
\cite{D:RNAM}). 

If $n \ge r-1$ the representation $\sigma$ in \eqref{7} is faithful
\cite{Brown}; thus it induces an isomorphism $\B_r^{(\epsilon n)}
\simeq \End_G(V^{\otimes r})$.

\subsection{Schur-Weyl duality in type $D$}\label{2.2a}
In type $D_{n/2}$ ($n$ even) the orthogonal group $\O_n$ is not
connected, and contains the connected semisimple group $\SO_n$
(special orthogonal group) as subgroup of index 2. In order to handle
this situation, Brauer (see also \cite{Grood}) defined a larger
algebra $\mathcal{D}_r^{(n)}$, spanned by the usual $r$-diagrams
previously defined, together with certain {\em partial} $r$-diagrams
on $2r$ vertices and $r-n$ edges, in which $n$ vertices in each of the
top and bottom rows are not incident to any edge, and showed that the
action of $\B_r^{(n)}$ can be extended to an action of this larger
algebra $\mathcal{D}_r^{(n)}$ on $V^{\otimes r}$. Thus we have
representations
\begin{equation}\label{so:1}
\begin{CD}
\C \SO_n @>\rho>> \End(V^{\otimes r}) @<\sigma<< \mathcal{D}_r^{(\epsilon n)}. 
\end{CD}
\end{equation}
Brauer showed that the actions of $\SO_n$
and $\mathcal{D}_r^{(n)}$ on $V^{\otimes r}$ satisfy Schur-Weyl
duality: 
\begin{align}
\rho(\C \SO_n) &= \End_{\mathcal{D}_r^{(\epsilon n)}}(V^{\otimes r}) 
\label{so:2}\\
\sigma(\mathcal{D}_r^{(\epsilon n)}) &= \End_{\SO_n}(V^{\otimes r}). 
\label{so:3}
\end{align}
The algebra in \eqref{so:2} is a {second} Schur algebra in type $D$, a
proper subalgebra of the algebra $\End_{\B_r^{(n)}}(V^{\otimes r})$
appearing in \eqref{7a} above.

\subsection{Generators and relations}
One can formulate the above statements of Schur-Weyl duality using
enveloping algebras, analogous to \ref{1.2}. This leads to a
presentation (see \cite{DGS}) of the symplectic and orthogonal Schur
algebras which is compatible with (a slight modification of) the usual
Serre presentation of the enveloping algebra $\U(\g)$, where $\g =
\sp_n$ ($n$ even) or $\so_n$.

\subsection{The quantum case}
There is a $q$-version of the Schur-Weyl duality considered in this
section, although not as developed as in type $A$. One needs to
replace the Brauer algebra by its $q$-analogue, the
Birman-Murakami-Wenzl (BMW) algebra (see \cite{BW}, \cite{Murakami}),
and replace the enveloping algebra by a suitable quantized enveloping
algebra.  One can think of the BMW algebra in terms of Kauffman's
tangle monoid; see \cite{Kauf}, \cite{HalvRam},
\cite{MortWass}. (Roughly speaking, tangles are replacements for
Brauer diagrams, in which one keeps track of under and over crossings,
subject to certain natural relations.)  There are applications of the
BMW algebra to knot theory, as one might imagine.

This leads to a $q$-analogue of the symplectic Schur algebras, in
particular, which have been studied by Oehms \cite{Oe}.

To the best of my knowledge, a $q$-analogue of the larger algebra
$\mathcal{D}_r^{(n)}$ ($n$ even) considered in \ref{2.2a} remains to
be formulated.

\section{The walled Brauer algebra}

\subsection{The algebra $\B_{r,s}^{(x)}$} \label{3.1}
This algebra was defined in 1994 in \cite{B}. It is the subalgebra of
$\B_{r+s}^{(x)}$ spanned by the set of $(r,s)$-diagrams. By
definition, an $(r,s)$-diagram is an $(r+s)$-diagram in which we
imagine a wall separating the first $r$ from the last $s$ columns of
vertices, such that:

 (a) all horizontal edges cross the wall;

 (b) no vertical edges cross the wall.

\noindent
An edge crosses the wall if its two vertices lie on opposite sides of
the wall. The multiplication in $\B_{r,s}^{(x)}$ is that of
$\B_{r+s}^{(x)}$.

Label the vertices on the top and bottom rows of an $(r,s)$-diagram by
the numbers $1, \dots, r$ to the left of the wall and $-1, \dots, -s$
to the right of the wall. Let $c_{i,-j}$ ($1 \le i \le r$; $1\le j \le
s$) be the diagram with a horizontal edge connecting vertices $i$ and
$-j$ on the top row and the same on the bottom row, and with all other
edges connecting vertex $k$ ($k \ne i, -j$) in the top and bottom
rows. It is easy to see that the walled Brauer algebra is generated by
the permutations it contains along with just one of the $c_{i,-j}$.
(Note that $c_{i,-j}$ is the $(r+s)$-diagram denoted by $c_{i,r+j}$ in
\ref{2.1}.)

\subsection{Dimension}
\begin{figure} [ht]
       \centering
       \includegraphics{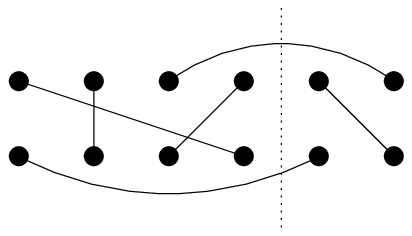}\qquad
       \includegraphics{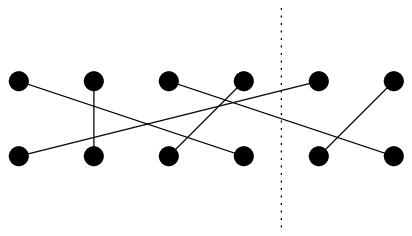}
       \caption{A $(4,2)$-diagram and its corresponding permutation,
       after applying {\em flip}.}
       \label{fig4}
       \end{figure}

What is the dimension of $\B_{r,s}^{(x)}$? One way to answer that
question is to consider the map, {\em flip}, from $(r+s)$-diagrams to
$(r+s)$-diagrams, defined by interchanging the top and bottom vertices
to the right of the imaginary wall.  For example, Figure \ref{fig4}
shows a $(4,2)$-diagram (to the left) and its corresponding
$6$-diagram, obtained from the left diagram by applying {\em flip}.
Note that {\em flip} is involutary: applying it twice gives the
original diagram back again.

One easily checks that the map {\em flip} carries $(r,s)$-diagrams
bijectively onto the set of $(r+s)$-diagrams with all edges
vertical. Such diagrams correspond with permutations of $r+s$ objects,
so the dimension of $\B_{r,s}^{(x)}$ is $(r+s)!$.

\subsection{Another view of $\B_{r,s}^{(x)}$}
The above correspondence between $(r,s)$-diagrams and permutations
gives another way to think of the multiplication in $\B_{r,s}^{(x)}$.
Given two $(r,s)$-diagrams $D_1$, $D_2$ let $D_1^\prime$, $D_2^\prime$
be their corresponding permutations obtained by applying {\em flip}.
Define a new (rather bizarre) composition on permutations as
follows. Given any two permutation diagrams $D_1^\prime$, $D_2^\prime$
(with $r+s$ columns of vertices) identify the first $r$ vertices of
the bottom row of $D_1^\prime$ with the first $r$ vertices of the top
row of $D_2^\prime$, and identify the last $s$ vertices of the top row
of $D_1^\prime$ with the last $s$ vertices of the bottom row of
$D_2^\prime$.  After removing loops and identified vertices this gives
a new permutation diagram $D_3^\prime$ in which the vertices in the
top (resp., bottom) row are the remaining top (bottom) row vertices
from the original diagrams.

Let $\delta(D_1^\prime, D_2^\prime)$ be the number of loops removed in
computing the composite permutation diagram $D_3^\prime$.  Define
multiplication of permutation diagrams by the rule
$$
 D_1^\prime \cdot D_2^\prime = x^{\delta(D_1^\prime, D_2^\prime)}
 D_3^\prime
$$
In other words, we are multiplying permutations by composing maps ``on
the right'' on one side of the wall, and ``on the left'' on the other
side (roughly speaking).  For example, Figure \ref{fig5} below shows
the computation of the composite diagram in the walled Brauer algebra
(left column) and the computation in terms of the corresponding
permutations (right column).  Figure \ref{fig6} shows the resulting
diagrams after the single loop and identified vertices have been
removed.
       \begin{figure} [ht]
       \centering
       \includegraphics{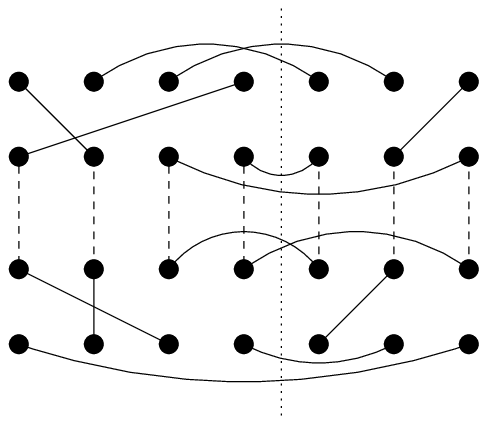}\qquad
       \includegraphics{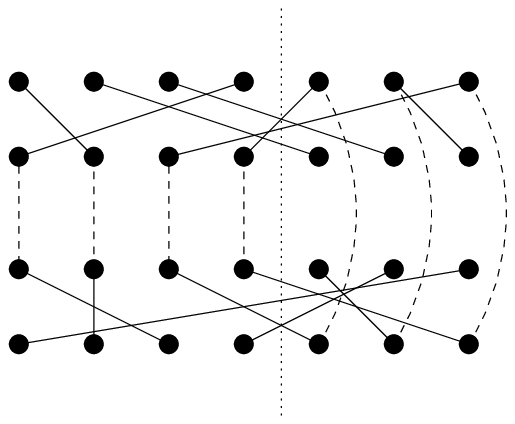}
       \caption{Composition of diagrams and permutations.
       The diagrams on the left correspond under {\em flip}
       with the permutations on the right.}
       \label{fig5}
       \end{figure}
       \begin{figure} [ht]
       \centering
       \includegraphics{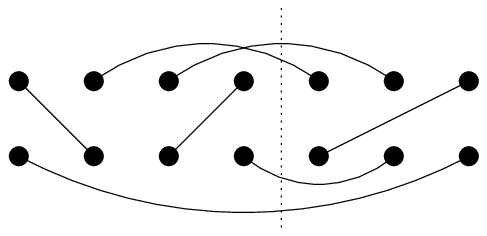}\qquad
       \includegraphics{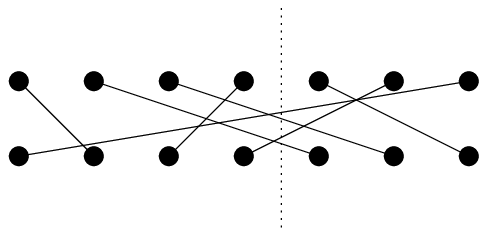}
       \caption{The corresponding diagrams resulting from Figure \ref{fig5}.
       The two diagrams correspond under {\em flip}.}
       \label{fig6}
       \end{figure}

One can check that $D_3^{\prime}$ coincides with $(D_1 \circ
D_2)^\prime$ and $\delta(D_1,D_2) = \delta(D_1^\prime, D_2^\prime)$.
In other words, {\em flip} defines an algebra isomorphism between the
algebra $\B_{r,s}^{(x)}$ and the algebra $\widetilde{\B}_{r,s}^{(x)}$
spanned by permutation diagrams with the multiplication defined above.
Note that in particular $\widetilde{\B}_{r,0}^{(x)} \simeq R[x]\S_r$ and  
$\widetilde{\B}_{0,s}^{(x)} \simeq (R[x]\S_s)^{\text{opp}}$.

\subsection{Schur-Weyl duality}
Consider mixed tensor space $V^{r,s} := V^{\otimes r} \otimes
V^{*\otimes s}$, where $V^*$ is the usual linear dual space of $V$.
Mixed tensor space is naturally a module for $\GL_n$, and one obtains
an action of $\B_{r,s}^{(n)}$ on $V^{r,s}$ simply by restricting the
action of $\B_{r+s}^{(n)}$, which acts the same on $V^{r,s}$ as it
does on $V^{\otimes(r+s)}$, since on restriction to $\O_n$ we have $V
\simeq V^*$.  Thus we have the following commutative diagram
\begin{equation}\label{wsw:1}
\begin{CD}
\C \GL_n @>\rho>> \End(V^{r,s}) @<\sigma<< \B_{r,s}^{(n)} \\
@A\iota AA & & @VV\iota' V \\
\C \O_n @>\rho>> \End(V^{\otimes(r+s)}) @<\sigma<< \B_{r+s}^{(n)}
\end{CD} 
\end{equation}
in which the vertical maps $\iota$, $\iota'$ are inclusion.  By
\cite{B}, the actions of $\GL_n$ and $\B_{r,s}^{(n)}$ on $V^{r,s}$
in the first row of the diagram satisfy Schur-Weyl duality:
\begin{align}
\rho(\C \GL_n) &= \End_{\B_{r,s}^{(n)}}(V^{r,s}) 
\label{wsw:2}\\
\sigma(\B_{r,s}^{(n)}) &= \End_{\GL_n}(V^{r,s}). 
\label{wsw:3}
\end{align}
The algebra in \eqref{wsw:2} is another Schur algebra $S(n; r,s)$ in
type $A$, studied in \cite{DiDo}.  These Schur algebras provide us
with a new family of quasihereditary algebras, generalizing the
classical Schur algebras, since $S(n; r,0) \simeq S(n,r)$.  In fact,
the $S(n; r,s)$ provide a new class of generalized Schur algebras in
the sense of Donkin \cite{Donkin:SA1}.  For fixed $n$, the family of
$S(n; r,s)$-modules as $r,s$ vary constitutes the family of all
rational representations of $\GL_n$.  Whence the name {\em rational}
Schur algebras for the $S(n; r,s)$.

When $n \ge r+s$, the representation $\sigma$ in the top row of
\eqref{wsw:1} above is faithful, so induces an isomorphism
$\B_{r,s}^{(n)} \simeq \End_{\GL_n}(V^{r,s})$.

\subsection{The quantum case}
Quantizations of the walled Brauer algebra have been defined and
studied in work of Halverson \cite{Halv}, Leduc \cite{Leduc},
Kosuda-Murakami \cite{KM}, and Kosuda \cite{Kosuda}.

\section{The deranged algebra}

\subsection{The problem}
One might wonder if there are versions of Schur-Weyl duality in which
the natural module $V$ is replaced by some other representation.
Perhaps the first choice would be to replace $V$ with the adjoint
module, {\em i.e.}, an algebraic group acting on its Lie algebra via
the adjoint representation.  The simplest instance of this would be
type $A$, where we consider the module $\sl_n^{\otimes r}$ as an
$\SL_n$-module, and ask for its centralizer algebra
$\End_{\SL_n}(\sl_n^{\otimes r}) = \End_{\GL_n}(\sl_n^{\otimes r})$.
(It makes no difference whether we regard $\sl_n$ as module for
$\SL_n$ or for $\GL_n$. We look at $\sl_n^{\otimes r}$ rather than
$\gl_n^{\otimes r}$ since $\gl_n$ is not simple as a $\GL_n$-module or
$\SL_n$-module.)

\subsection{Relation with Brauer algebras}
Even though this is a question about type $A$, its solution is
intimately connected with the walled Brauer algebra.  Here is a brief
outline of the solution to this problem, recently obtained in
\cite{BD}. The main idea is to utilize the decomposition $\gl_n =
\sl_n \oplus \C$ (as $\SL_n$ or $\GL_n$ module) to write
\begin{equation}
\gl_n^{\otimes r} = \sl_n^{\otimes r} \oplus \textstyle \bigoplus_{0
\le t < r} \sl_n^{\otimes t}
\end{equation}
where the left-hand-side identifies with $V^{r,r} = V^{\otimes r}
\otimes V^{*\otimes r}$ via the natural isomorphism $V\otimes V^*
\simeq \End(V) \simeq \gl_n$. Thus it follows that our desired tensor
space $\sl_n^{\otimes r}$ is isomorphic with a direct summand (as a
$\GL_n$ or $\SL_n$-module) of the mixed tensor space $V^{r,r}$, and
for $n\ge 2r$ the centralizing algebra $\mathcal{C} =
\End_{\GL_n}(\sl_n^{\otimes r})$ is obtainable as a certain subalgebra
$e \End_{\GL_n}(V^{r,r}) e$ where $e$ is the idempotent corresponding
to the projection onto the summand $\sl_n^{\otimes r}$.  Thus
\begin{equation}\label{9}
\mathcal{C} \simeq e\B_{r,r}^{(n)}e \qquad(n\ge 2r)
\end{equation}
where $e = \prod_{i}(1 - n^{-1} c_{i,-i})$ (notation of \ref{3.1}).
The algebra in \eqref{9} has a basis consisting of all elements of the
form $eDe$ as $D$ ranges over the set of $(r,r)$-diagrams with no
horizontal edges connecting $i$ to $-i$. Such diagrams correspond
under {\em flip}, after inverting the signs labeling the bottom row,
with {\em derangements} of $2r$ objects.  (A derangement is a
permutation having no fixed points.)  Let $N(k) = $ the number of
derangements of $k$ objects, then we have by the Inclusion-Exclusion
Principle (see \cite[2.2.1 Example] {Stanley})
\begin{equation}
N(k) = \sum_{j=0}^k (-1)^{k-j} \binom{k}{j} \, j!
\end{equation}
which is the nearest integer to $k!/\mathsf{e}$ ($\mathsf{e} =
2.7182818\dots$).  Thus the dimension of the centralizer algebra
$e\B_{r,r}^{(n)}e$ is given by $N({2r})$. Because of this connection
with derangements, the algebra $e\B_{r,r}^{(n)}e$ is known as the {\em
deranged algebra}.

In \cite{BD} explicit formulas are obtained for the number of times a
given simple module $L$ appears as a summand of $\sl_n^{\otimes r}$.
In particular, it is shown that
$$
N(r) = \text{the multiplicity of the trivial module in $\sl_n^{\otimes r}$}
$$
and 
$$
N(r-1) = \text{the multiplicity of $\sl_n$ in $\sl_n^{\otimes r}$}
$$
demonstrating that the combinatorics of derangement numbers is
inherent in this theory.

\subsection{Schur-Weyl duality}
At least when $n \ge 2r$, the actions of $\GL_n$ and $\mathcal{C}$
on $\sl_n^{\otimes r}$ satisfy Schur-Weyl duality:
\begin{align}
\rho(\C \GL_n) &= \End_{e\B_{r,r}^{(n)}e}(\sl_n^{\otimes r}) 
\label{dsw:2}\\
\sigma(e\B_{r,r}^{(n)}e) &= \End_{\GL_n}(\sl_n^{\otimes r}).
\label{dsw:3}
\end{align}
This will almost certainly hold for all $n$, $r$.

\end{document}